\documentclass[3p,preprint]{elsarticle}

\usepackage{amssymb}

\usepackage{amsmath} 

\usepackage{lineno}

\usepackage{amssymb,amsthm,MnSymbol}
\usepackage{color}

\usepackage{algorithm}
\usepackage{algorithmic}

\usepackage{subfig}

\usepackage{amscd}
\usepackage{mathrsfs}
\usepackage[colorlinks]{hyperref}

\usepackage{enumitem}

\definecolor{darkgreen}{rgb}{0.0,0.5,0.0}               
\definecolor{lightblue}{rgb}{0.3296, 0.6648, 0.8644}    
\definecolor{shadowcolor}{rgb}{0.0000, 0.0000, 0.6179}  
\definecolor{bulletcolor}{rgb}{ 0.8441, 0.1582, 0.0000} 
\definecolor{deepskyblue}{rgb}{0, 0.74902, 1.0}
\definecolor{royalblue}{rgb}{0.254901960784314,   0.411764705882353,   0.882352941176471}
\definecolor{dodgerblue}{rgb}{0.11765, 0.56471, 1.0}
\definecolor{bordercolor}{rgb}{0,0,.2380}        
\definecolor{internationalorange}{rgb}{1.0, 0.31, 0.0}

\newcommand{\ansA}[1]{{\color{black}#1}}
\newcommand{\ansB}[1]{{\color{black}#1}}
\usepackage[normalem]{ulem}

\renewcommand\sout{\bgroup\protect\markoverwith{\textcolor{black}{\protect\rule[0.5ex]{2pt}{0.8pt}}}\protect\ULon}

\biboptions{authoryear}

\usepackage[figuresright]{rotating}

\usepackage{epstopdf}

\begin{document}
\begin{frontmatter}

\author[addr1]{I. Shevchenko\corref{cor1}}
\ead{i.shevchenko@imperial.ac.uk}
\author[addr1,addr2]{P. Berloff}
 
\cortext[cor1]{Corresponding author at:}
\address[addr1]{Department of Mathematics, Imperial College London, Huxley Building, 180 Queen's Gate, London, SW7 2AZ, UK}
\address[addr2]{Institute of Numerical Mathematics of the Russian Academy of Sciences, Moscow, Russia}


\title{A method for preserving nominally-resolved flow patterns in low-resolution ocean simulations}

\begin{abstract}
Accurate representation of large-scale flow patterns in low-resolution ocean simulations is one of the most challenging problems 
in ocean modelling.
The main difficulty is to correctly reproduce effects of unresolved small scales on the resolved large scales.
For this purpose, most of current research is focused on development of parameterisations accounting for the small scales.
In this work we propose an alternative to the mainstream ideas by showing how to reconstruct a dynamical system 
from the available reference solution data (our proxy for observations) and, then, how to use this system for modelling not only large-scale but 
also nominally-resolved flow patterns at low resolutions.
The proposed methodology robustly allows to retrieve a system of equations underlying reduced dynamics of the observed data.
Remarkably, its solutions simulate not only
large-scale but also small-scale flow features, which can be nominally resolved by the low-resolution grid.
\end{abstract}


\end{frontmatter}

\section{Introduction}
It is typical of low-resolution ocean simulations to have significantly distorted or even absent large-scale flow structures 
that are otherwise present in the high-resolution simulations.
This failure is due to missing effects of the small scales, which are not adequately resolved in low-resolution simulations.
To mitigate the problem, many parameterisations for both comprehensive and idealized ocean models have been proposed (e.g.,~\citet{GentMcwilliams1990,DuanNadiga2007,Frederiksen_et_al2012,JansenHeld2014,
PortaMana_Zanna2014,CooperZanna2015,
Grooms_et_al2015,Berloff_2015,Berloff_2016,Berloff_2018,Danilov_etal_2019,Ryzhov_etal_2019,Juricke_etal_2020a,Juricke_etal_2020b,
CCHWS2019_1,Ryzhov_etal_2020,
CCHWS2019_3,CCHWS2020_4,CCHPS2020_J2}), but overall the problem remains largely unresolved for several reasons.
First, defining the small and large scales is ambiguous, because they are not separated by a clear spectral gap or otherwise.
Second, definition of the small and large scales \ansA{should be} consistent with the specific resolving capabilities 
of a low-resolution model in which their interactions are to be parameterized;
a new definition of eddies as field error of the employed model has been proposed in~\citep{BRS2020_J1}.
Third, dynamical interactions across the scales are remarkably complex, as well as spatially inhomogeneous and non-stationary. 

In this paper the problem is approached from a different direction: instead of parameterizing small-scale effects, we retrieve an underlying dynamical system and 
use it to model evolution of the nominally-resolved flow patterns 
(the ones that can be properly resolved on the coarse grid) at low resolutions.
Although the basic idea has long research history, our application of it is novel, and the proposed methodology has many novel features.
Let us first discuss below the relevant background.

Retrieving reduced equations underlying the observed flow evolution is one of the most challenging problems in dynamical systems (e.g.,\citep{AguirreLetellier_2009,Brunton_et_al_2016}).
Although, this field has been researched for decades, most of the efforts used low-dimensional dynamical systems with 3-5 degrees of freedom 
(e.g., ~\citet{Brunton_et_al_2017,MangiarottiHuc_2019}), and even this turned out difficult. 
This is, because with larger number of degrees of freedom, there are so many types of models and various subtleties that investigating them becomes a big task. 
The other problem is about frequent numerical instabilities of the retrieved dynamical systems.
This implies that applying known methodologies for thousands of degrees of freedom, typical for describing low-resolution oceanic flows, is unfeasible. 

For developing and testing the approach, we considered an intermediate-complexity, quasigeostrophic, eddy-resolving model of the wind-driven midlatitude 
ocean circulation --- this is a respected and widely used 
(e.g.,~\citet{SiegelEtAl2001,Karabasov_et_al2009,SB_FLUIDS_2016,SBGR_JFM_2016} and references there in) 
paradigm for process studies involving large-/small-scale turbulent interactions and their parameterizations.
To mitigate the model size problem, we applied the Empirical Orthogonal Function (EOF) analysis ~\citep{Preisendorfer1988,Hannachi_et_al2007} 
to the reference flow, defined here as the high-resolution solution subsampled on a coarse grid, and reconstructed a dynamical 
system for the evolving Principal Components (PCs) corresponding to the leading EOFs.
\ansB{Successful examples of reduced-order modelling with EOF-PC description can be found in~\citep{KB_2015,KCB_2018}.}
\ansB{Other types of space reduction are possible and can improve the outcome even further, but they are not
considered in this study.}
To resolve the problem with numerical instabilities, we used adaptive nudging 
methodology, which is an upgraded extension of the nudging method proposed in ~\citep{ShevchenkoBerloff_2021}. 

\section{The method\label{sec:method}}

The main objective of this study is to reconstruct a dynamical system from the reference solution (say $\mathbf{x}(t)$, $\mathbf{x}\in\mathbb{R}^n$).
This dynamical system is meant to correctly simulate the reference large-scale flow patterns on a low-resolution grid. 
The full dimensionality of the problem is the total number of the grid nodes, and 
for the reconstructed dynamical system we aim to reduce it by orders of magnitude via EOFs/PCs decomposition of the reference data.
Next, we formulate general dynamical system in terms of the leading PCs:
\begin{equation}
\mathbf{y}'(t)=\mathbf{F}(\mathbf{y}),\quad \mathbf{y}\in\mathbb{R}^m,\quad t\in[0,\widetilde{T}],\quad m<<n \, ,
\label{eq:ode1} 
\end{equation}
where the PCs are combined in the vector and denoted by $\mathbf{y}(t)$.
In our case the dimensionality has been eventually reduced by three orders of magnitude (from $n=16441$ to $m=30$).

In ~\eqref{eq:ode1} we used 30 leading PCs that captured 98\% of the reference flow variance.
The right hand side of ~\eqref{eq:ode1} is approximated with polynomial of order two in all the variables, $\mathbf{P}(\mathbf{y})$, and with the Fourier series, $\mathcal{F}(\mathbf{y})$, containing 50 leading harmonics:
\begin{equation}
\mathbf{F}(\mathbf{y})\approx \mathbf{P}(\mathbf{y}) + \mathcal{F}(\mathbf{y}) \, .
\label{eq:poly2} 
\end{equation}
Note that $\mathbf{F}(\mathbf{y})$ can be approximated differently, and its optimal choice (beyond the scope of this work) is a challenge for the dynamical system reconstruction.
Without proper information for tailoring the right hand side more specifically, a polynomial expansion is justified by the Weierstrass approximation theorem, 
while the use of the Fourier series allows one to approximate the mean flow more accurately.
We will get back to this choice when discussing the results.

In order to find the coefficients in $\mathbf{P}(\mathbf{y})$ and Fourier series, we use the method of least squares. 
Having approximated $\mathbf{F}(\mathbf{y})$ up to a given order of accuracy, one can solve the reconstructed dynamical system
\begin{equation}
\mathbf{z}'(t)=\mathbf{P}(\mathbf{z}) + \mathcal{F}(\mathbf{z}),\quad \mathbf{z}\in\mathbb{R}^m,\quad t\in[0,T],\quad T>\widetilde{T} \, . 
\label{eq:ode3} 
\end{equation}
Note that this system is integrated over a time interval which is longer (here, 2 times longer) than that of the original system~\eqref{eq:ode1}, since the purpose of the reconstructed system is to reproduce the original flow dynamics well beyond the known data record.
In all further simulations we will have $\widetilde{T}=$2 years and $T=$4 years.
However, an accurate approximation of $\mathbf{F}(\mathbf{y})$ does not guarantee that system~\eqref{eq:ode3} can be easily solved, 
because the integration errors can quickly contaminate the solution and result in severe numerical instability --- this is what actually happened in our case.
In order to stabilize the numerical integration, we used the nudging methodology ~\citep{ShevchenkoBerloff_2021}:
\begin{equation}
\mathbf{z}'(t)=\mathbf{P}(\mathbf{z}) + \mathcal{F}(\mathbf{z})+
\eta\left(\frac{1}{N}\sum\limits_{k\in\mathcal{U}(\mathbf{z}(t))}\mathbf{y}(t_k)-\mathbf{z}(t)\right),\quad t\in[0,T] \, ,
\label{eq:ode3_nudging}
\end{equation}
where $\mathcal{U}(\mathbf{z}(t))$ is a neighbourhood of $\mathbf{z}(t)$,
and index $k$ is the timestep of the corresponding PC $\mathbf{y}(t_k)$; the timestep of the PC is the timestep with which the actual data for the 
EOF analysis was generated.
The neighbourhood is computed \ansA{in $l_2$ norm} as the average of $N=5$ points nearest to the solution $\mathbf{z}(t)$.
\ansA{Note that the number of neighbourhood points is a parameter, and its sensitivity should be explored and taken into account 
for each application of the proposed methodology.}

Having solved equation~\eqref{eq:ode3_nudging}, we approximated the reference solution by using the leading EOF-PC pairs as follows:
\begin{equation}
\mathbf{x}(t)\approx\sum\limits^m_{i=1}z_i(t)\mathbf{E}_i \, ,
\label{eq:eof_pc}
\end{equation}
with $\mathbf{E}_i$ and $z_i$ being the $i$-th EOF and PC, respectively.

Note that $N$ in equation~\eqref{eq:ode3_nudging} can be made time-dependent and adaptive, like the nudging coefficient $\eta$, which is the other important parameter.
In order to make the numerical integration stable with the Euler method, we used the following adaptive nudging:   
\begin{equation}
\eta(t_i)=\left\{ \begin{array}{ll}
\displaystyle
\eta(t_{i-1})+\eta_h & \text{if } \sigma(\mathbf{z}(t_i))>\max\nolimits_{t\in[0,\widetilde{T}]}\sigma(\mathbf{y}(t)), \\
\eta(t_{i-1})-\eta_h & \text{if } \sigma(\mathbf{z}(t_i))\le\max\nolimits_{t\in[0,\widetilde{T}]}\sigma(\mathbf{y}(t)),\quad i=1,2,\ldots \\
0 & \text{if } \eta(t_{i-1})-\eta_h<0.
\end{array}\right.
\label{eq:eta}
\end{equation}
with $\sigma$ being the standard deviation, $\eta_h=0.001$, and $\eta(t_0)=0$. 

We opted out for an adaptive nudging, as it keeps the system within a neighbourhood of the phase space region occupied by the reference solution.
As an alternative, a constant $\eta$ can be also used with some tuning and caution, keeping in mind that its small value may not be enough for keeping 
the solution within the right region and its large value may result in an over-stabilized solution with suppressed flow variability (slow flow dynamics).

\section{Multilayer quasi-geostrophic model\label{sec:qg}}

We consider a 3-layer quasi-geostrophic (QG) model with forcing and dissipation for the evolution of the potential vorticity (PV) anomaly $\mathbf{q}=(q_1,q_2,q_3)$ in domain $\Omega$~\citep{Pedlosky1987}:
\begin{equation}
\partial_t q_j+\mathrm{J}(\psi_j,q_j+\beta y)=\delta_{1j}F_{\rm w}-\delta_{j3}\,\mu\nabla^2\psi_j+\nu\nabla^4\psi_j,\quad j=1,2,3\, ,
\label{eq:pve}
\end{equation}
where $\mathrm{J}(f,g)=f_xg_y-f_yg_x$, $\delta_{ij}$ is the Kronecker symbol, and $\boldsymbol{\psi}=(\psi_1,\psi_2,\psi_3)$ is the velocity streamfunction in three layers.
The planetary vorticity gradient is $\beta=2\times10^{-11}\, {\rm m^{-1}\, s^{-1}}$,
the bottom friction parameter is $\mu=4\times10^{-8}\, {\rm s^{-1}}$, and the lateral eddy viscosity is $\nu=50\, {\rm m^2\, s^{-1}}$.
The asymmetric wind curl forcing, driving the double-gyre ocean circulation, is given by
\[
\displaystyle
F_{\rm w}=\left\{ \begin{array}{ll}
\displaystyle
-1.80\,\pi\,\tau_0\sin\left(\pi y/y_0\right), & y\in[0,y_0), \\
{\color{white}-}2.22\,\pi\,\tau_0\sin\left(\pi (y-y_0)/(L-y_0)\right), & y\in[y_0,L],\\
\end{array}\right.
\]
with the wind stress amplitude $\tau_0=0.03\, {\rm N\, m^{-2}}$ and the tilted zero forcing line
$y_0=0.4L+0.2x$, $x\in[0,L]$.
The computational domain $\Omega=[0,L]\times[0,L]\times[0,H]$ is a closed, flat-bottom basin with $L=3840\, \rm km$, and the total depth
$H=H_1+H_2+H_3$ given by the isopycnal fluid layers of depths (top to bottom): 
$H_1=0.25\, \rm km$, $H_2=0.75\, \rm km$, $H_3=3.0\, \rm km$. 

The PV anomaly $\boldsymbol{q}$ and the velocity streamfunction $\boldsymbol{\psi}$ are coupled through the system of elliptic equations:
\begin{equation}
\boldsymbol{q}=\nabla^2\boldsymbol{\psi}-{\bf S}\boldsymbol{\psi} \, ,
\label{eq:pv}
\end{equation}
with the stratification matrix
\[
{\bf S}=\left(\begin{array}{lll}
   {\color{white}-}1.19\cdot10^{-3} & -1.19\cdot10^{-3}                  &  {\color{white}-}0.0                  \\
  -3.95\cdot10^{-4}                 &  {\color{white}-}1.14\cdot10^{-3}  & -7.47\cdot10^{-4}                  \\
          {\color{white}-}0.0          & -1.87\cdot10^{-4}                  &  {\color{white}-}1.87\cdot 10^{-4} \\
\end{array}\right).
\]
\noindent
The stratification parameters are given in units of $\rm km^{-2}$ and chosen so, that the first and second Rossby deformation 
radii are $Rd_1=40\, {\rm km}$ and $Rd_2=23\, {\rm km}$, respectively; \ansA{the choice of these parameters is typical for the North Atlantic,
as it allows to simulate a more realistic but yet idealized eastward jet extension of the western boundary currents}.

System~(\ref{eq:pve})-(\ref{eq:pv}) is augmented with the integral mass conservation constraint~\citep{McWilliams1977}:
\begin{equation}
\partial_t\iint\limits_{\Omega}(\psi_j-\psi_{j+1})\ dydx=0,\quad j=1,2
\label{eq:masscon}
\end{equation}
with the zero initial condition, and with the partial-slip lateral boundary condition~\citep{HMG_1992}:
\begin{equation}
\left(\partial_{\bf nn}\boldsymbol{\psi}-\alpha^{-1}\partial_{\bf n}\boldsymbol{\psi}\right)\Big|_{\partial\Omega}=0 \, ,
\label{eq:bc}
\end{equation}
where $\alpha=120\, {\rm km}$ is the partial-slip parameter, and $\bf n$ is the normal-to-wall unit vector;
\ansA{no-flow-through boundary condition is also implemented (as part of the elliptic solver)}.
\ansA{The value of the parameter $\alpha$ is chosen based on the study by~\citet{SB2015}, where it has been shown that smaller values of $\alpha$
inhibit the eastward jet extension penetration length and volume transport, while larger values have much less pronounced influence on the jet.
As with other governing parameters used in this study, our choice of $\alpha$ is justified by a more realistic eastward jet.}

\begin{figure}[H]
\centering
\begin{tabular}{ccccc}
& \hspace*{0.5cm}\begin{minipage}{0.1\textwidth} {\bf (a)} \end{minipage} & 
\hspace*{1cm}\begin{minipage}{0.1\textwidth} {\bf (b)} \end{minipage} &
\hspace*{1cm}\begin{minipage}{0.1\textwidth} {\bf (c)} \end{minipage} &
\hspace*{1cm}\begin{minipage}{0.1\textwidth} {\bf (d)} \end{minipage}\\
\hspace*{-1.5cm}\begin{minipage}{0.02\textwidth}\rotatebox{90}{$t=2$ years}\end{minipage}  &
\hspace*{-1cm}\begin{minipage}{0.24\textwidth}\includegraphics[scale=0.1]{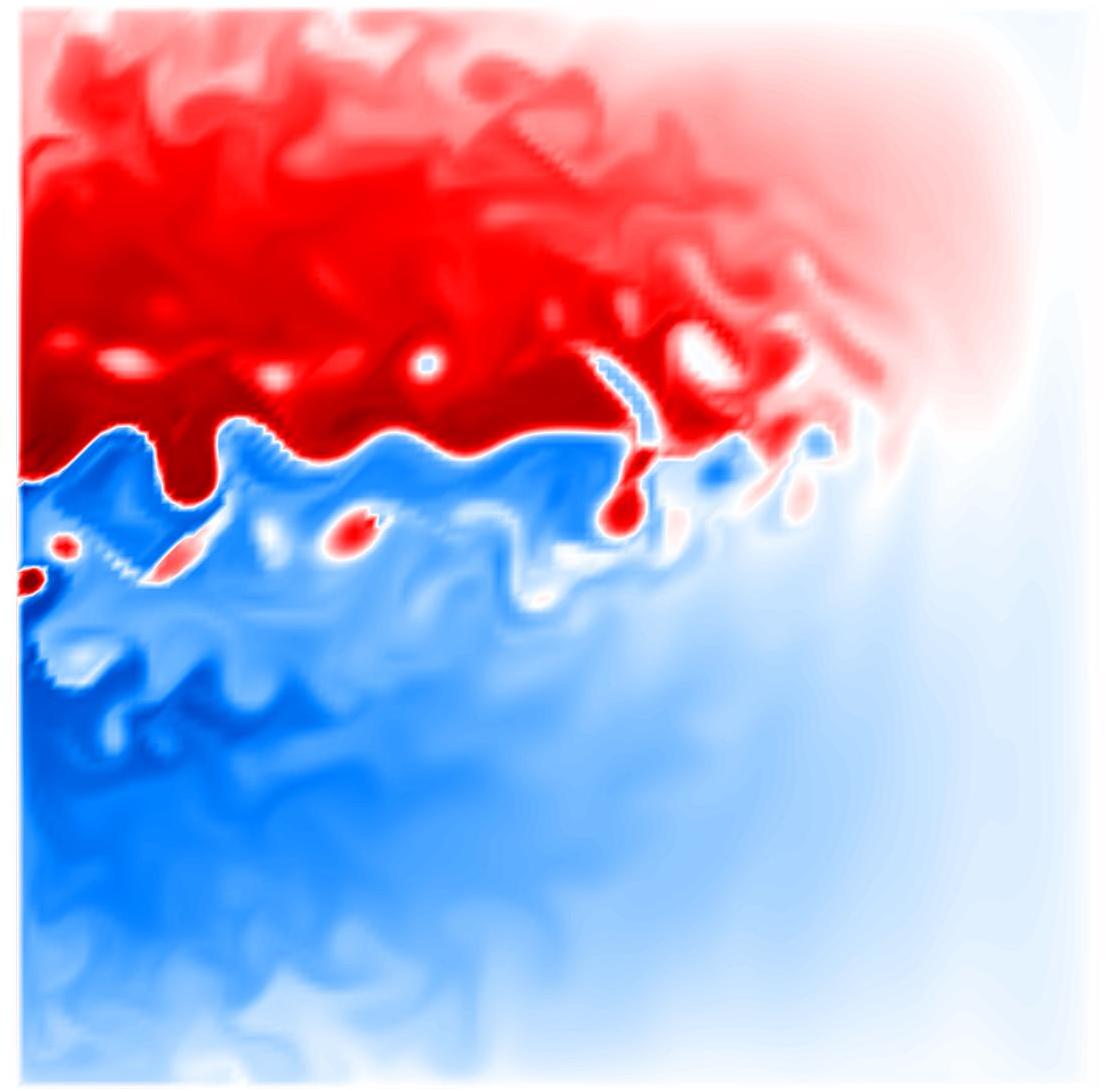}\end{minipage} &
\hspace*{0cm}\begin{minipage}{0.24\textwidth}\includegraphics[scale=0.1]{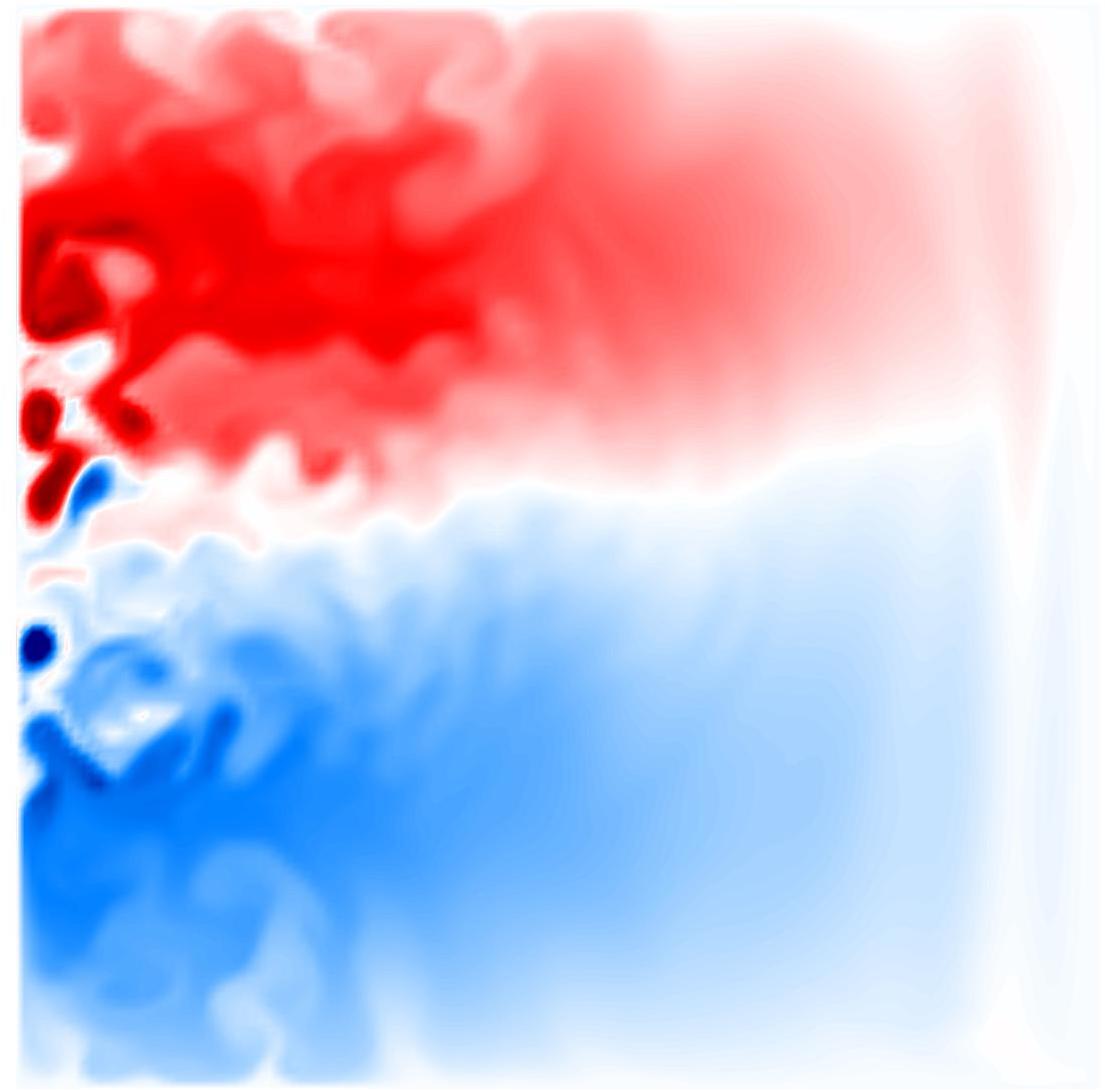}\end{minipage} &
\hspace*{0cm}\begin{minipage}{0.24\textwidth}\includegraphics[scale=0.1]{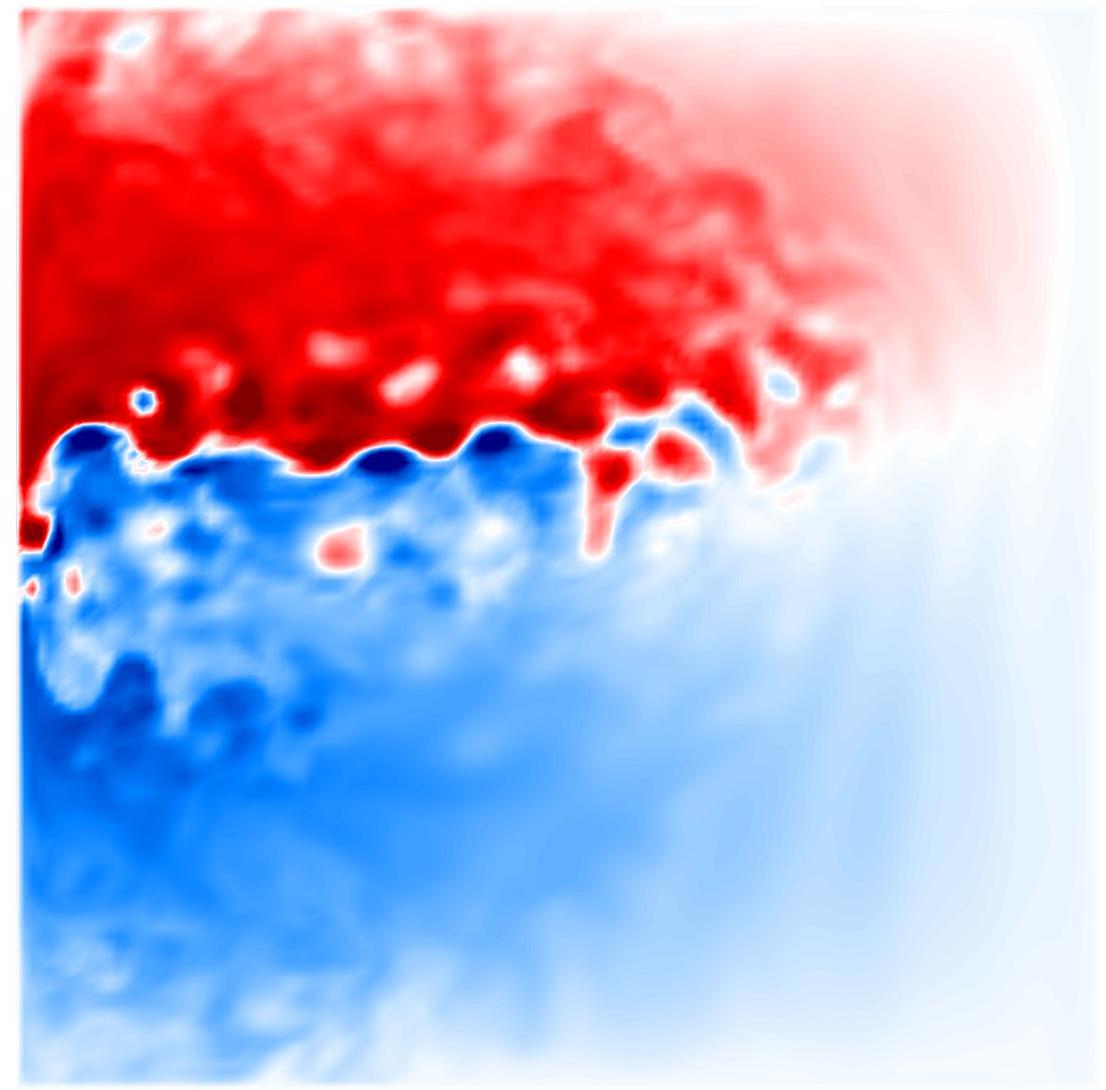}\end{minipage} &
\hspace*{0cm}\begin{minipage}{0.24\textwidth}\includegraphics[scale=0.1]{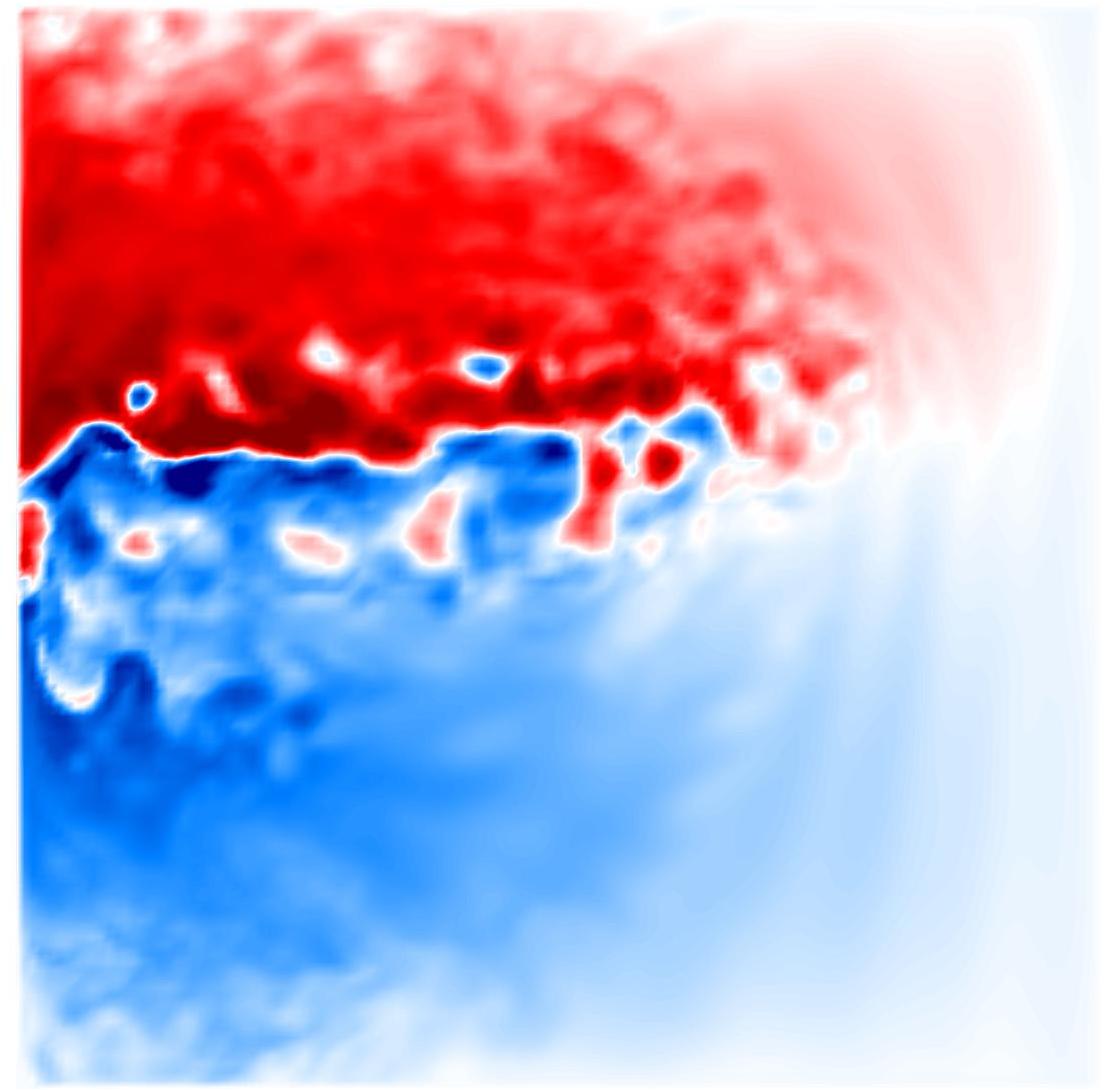}\end{minipage}\\
& & & & \\[-0.35cm]
\hspace*{-1.5cm}\begin{minipage}{0.02\textwidth}\rotatebox{90}{$t=4$ years}\end{minipage}  &
\hspace*{-1cm}\begin{minipage}{0.24\textwidth}\includegraphics[scale=0.1]{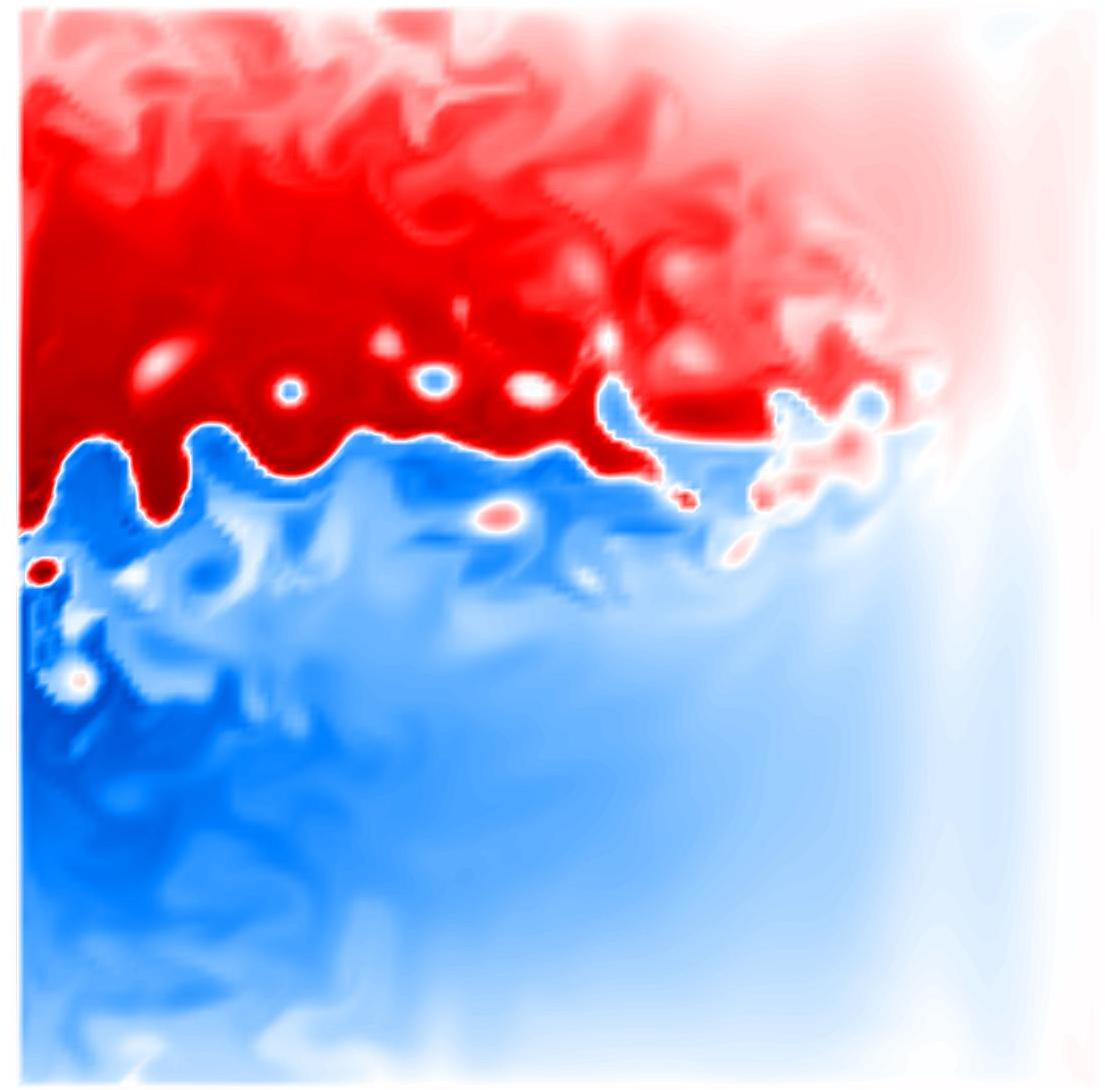}\end{minipage} &
\hspace*{0cm}\begin{minipage}{0.24\textwidth}\includegraphics[scale=0.1]{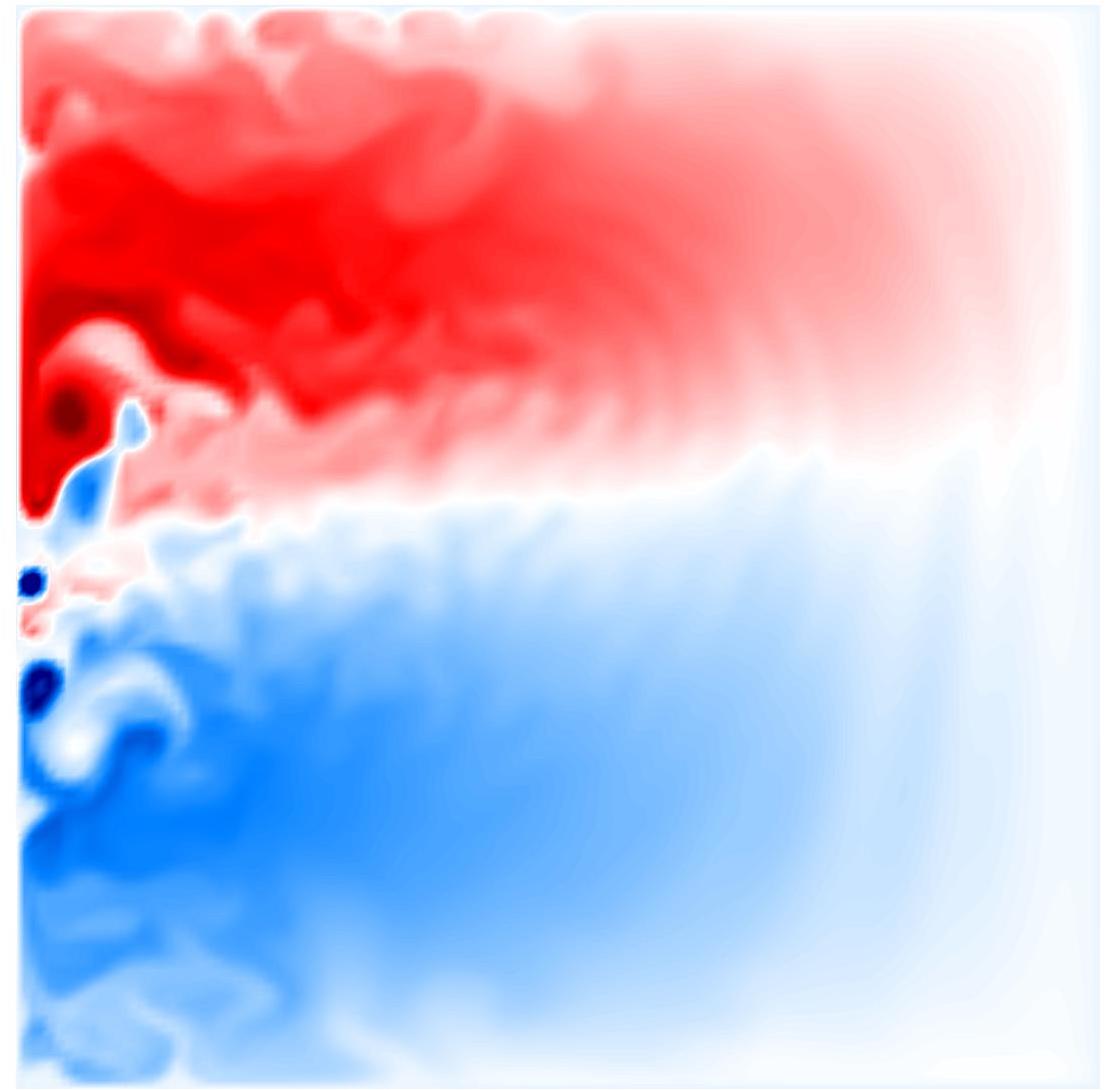}\end{minipage} &
\hspace*{0cm}\begin{minipage}{0.24\textwidth}\includegraphics[scale=0.1]{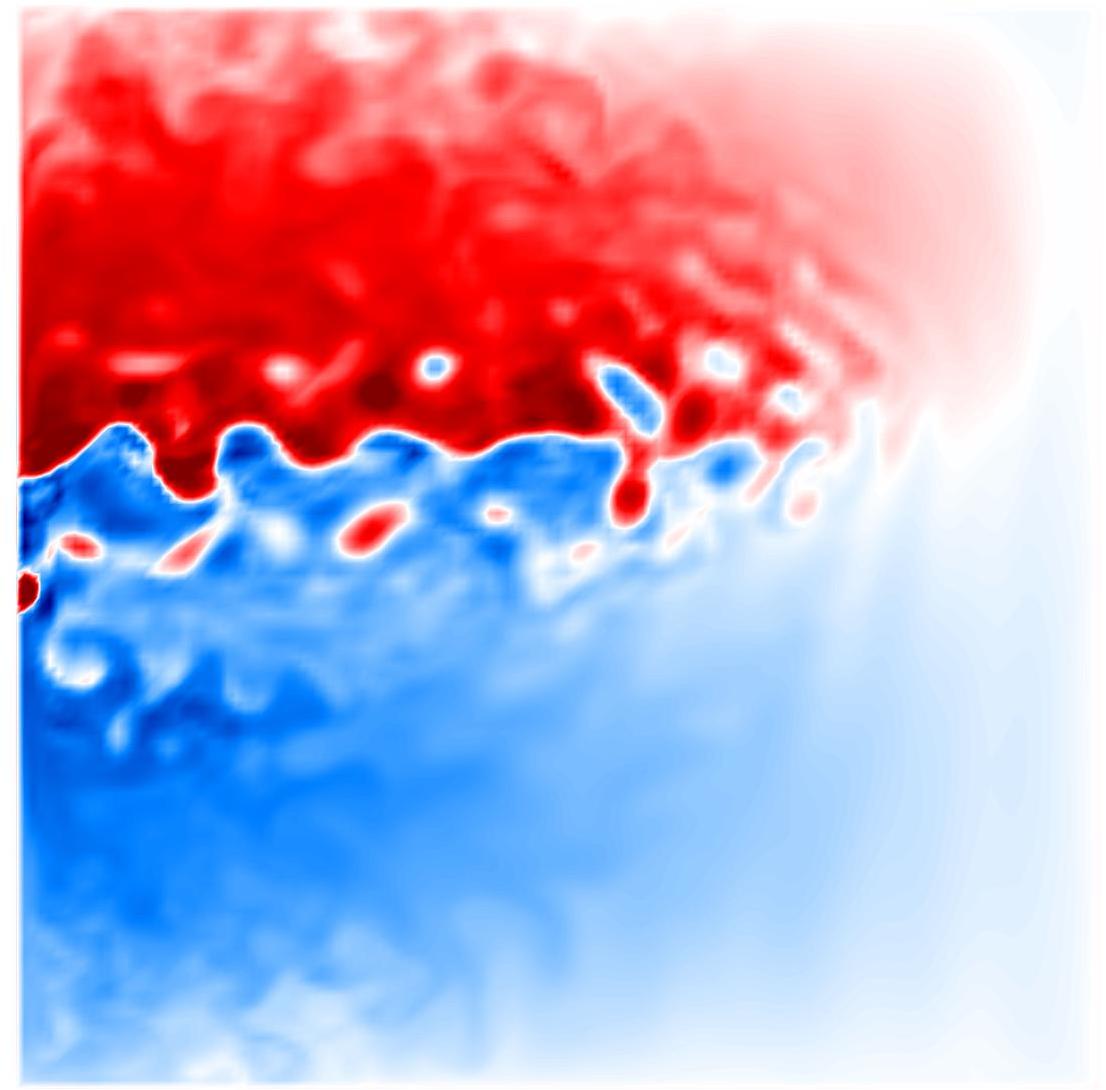}\end{minipage} &
\hspace*{0cm}\begin{minipage}{0.24\textwidth}\includegraphics[scale=0.1]{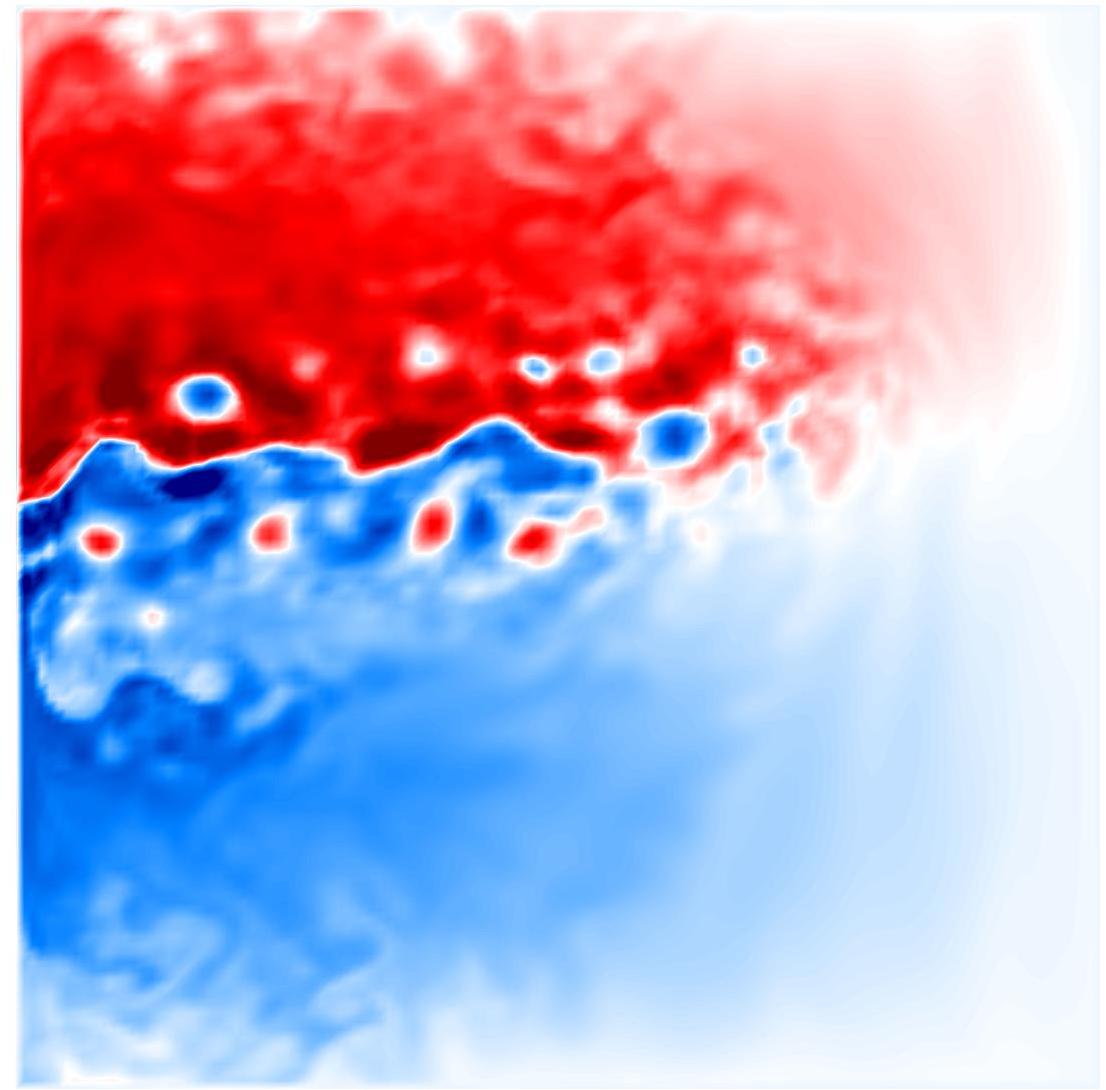}\end{minipage}\\
& & & & \\[-0.35cm]
\hspace*{-1.5cm}\begin{minipage}{0.02\textwidth}\rotatebox{90}{4-year average}\end{minipage}  &
\hspace*{-1cm}\begin{minipage}{0.24\textwidth}\includegraphics[scale=0.1]{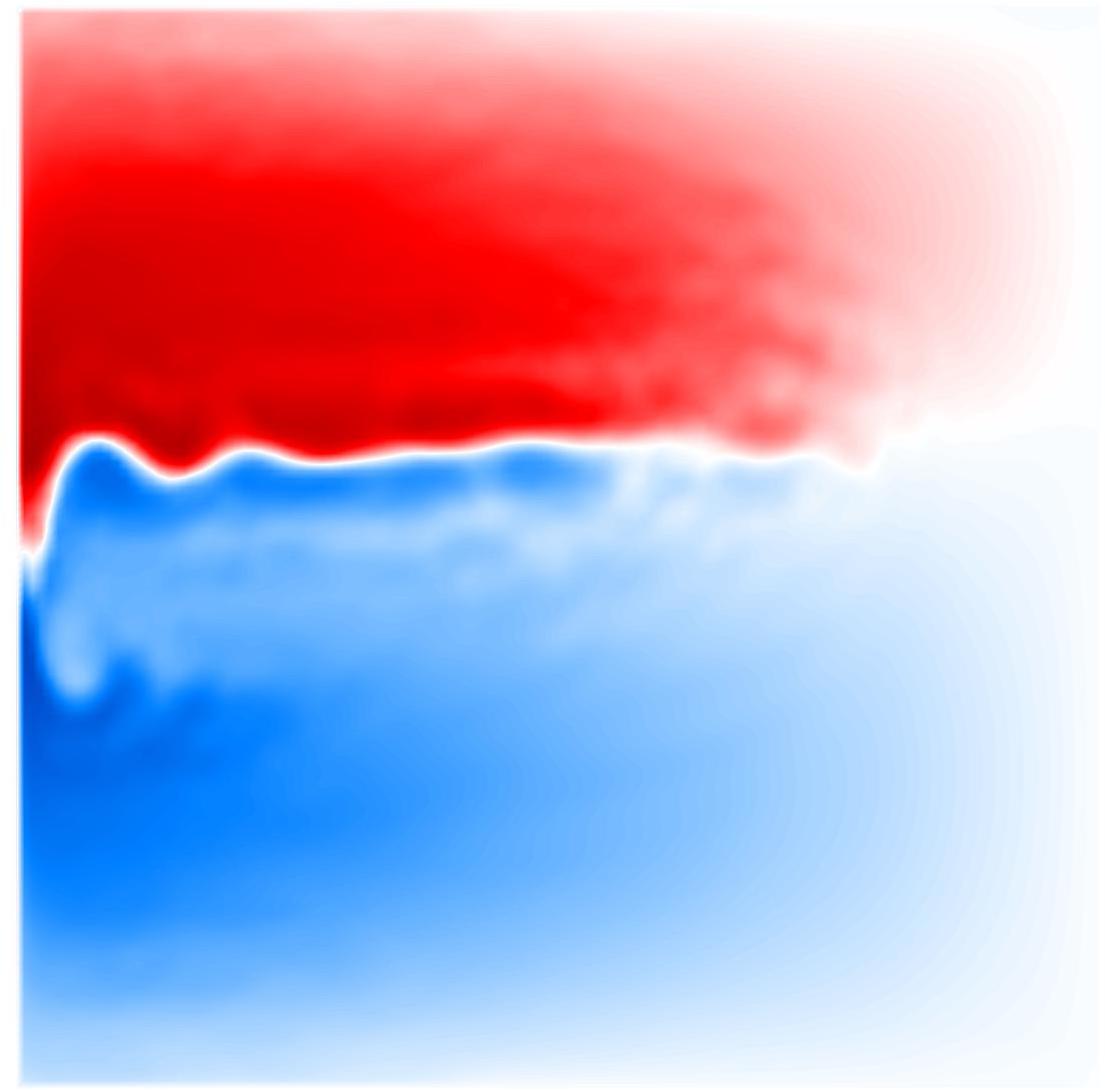}\end{minipage} &
\hspace*{0cm}\begin{minipage}{0.24\textwidth}\includegraphics[scale=0.1]{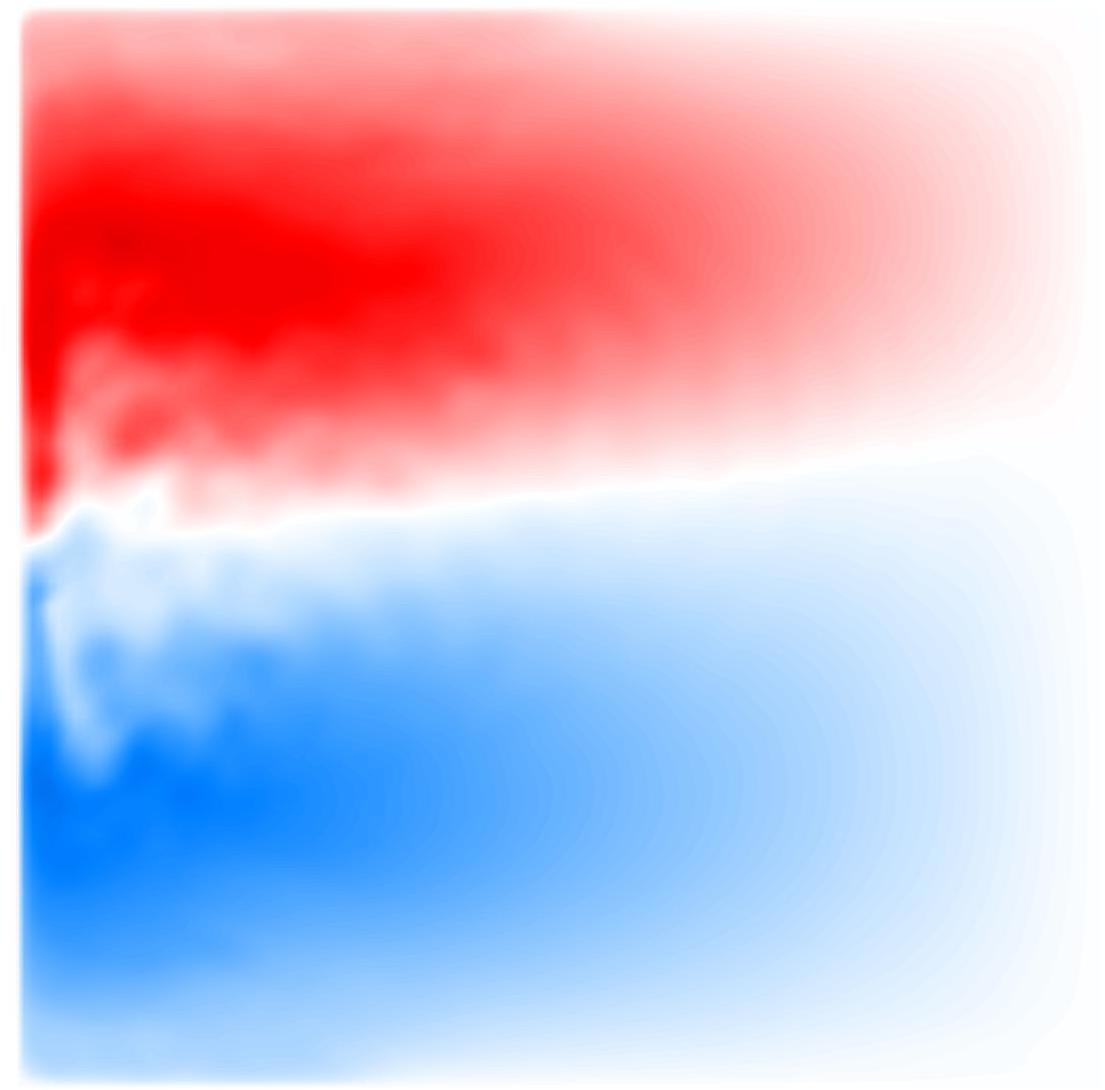}\end{minipage} &
\hspace*{0cm}\begin{minipage}{0.24\textwidth}\includegraphics[scale=0.1]{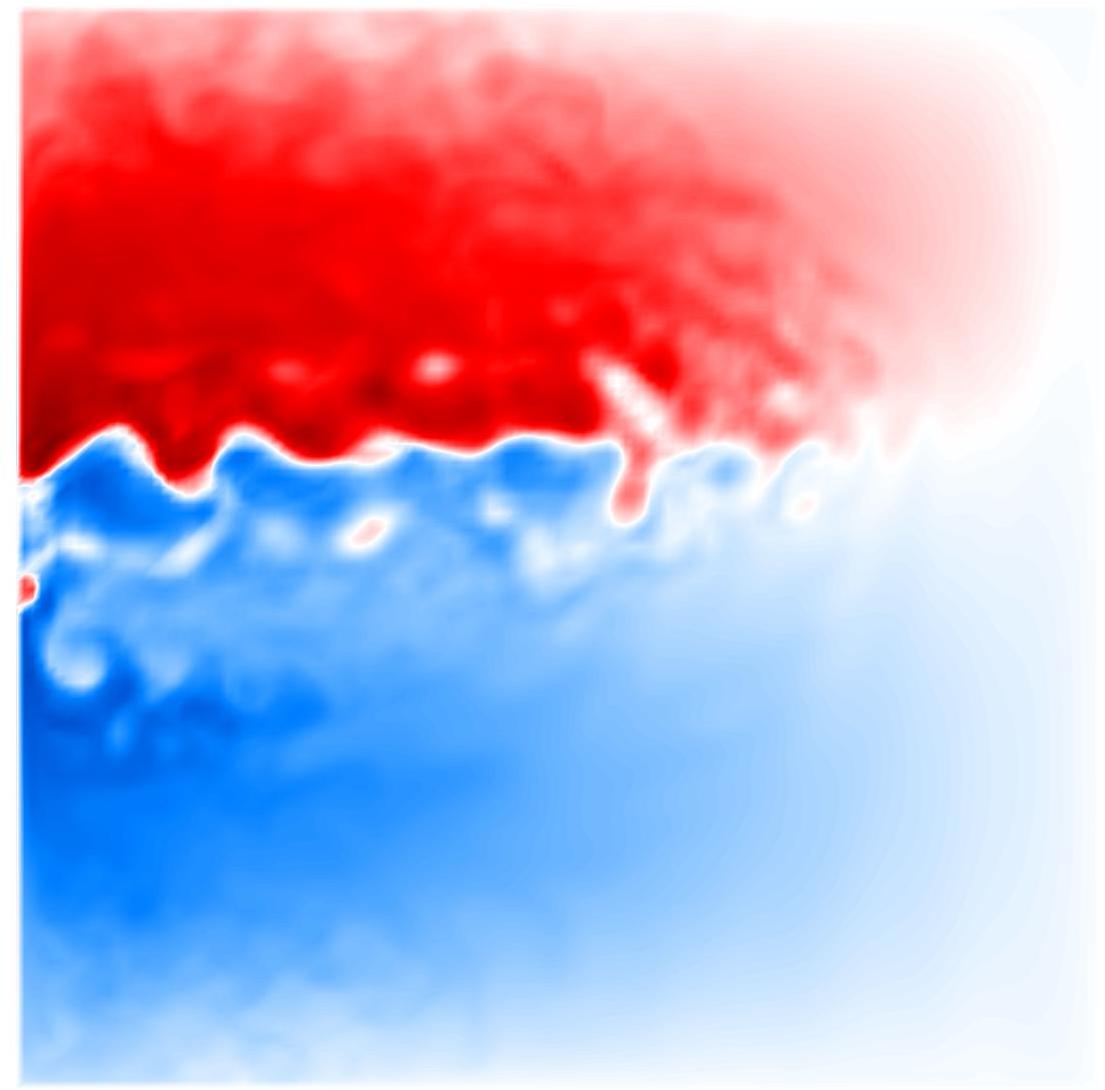}\end{minipage} &
\hspace*{0cm}\begin{minipage}{0.24\textwidth}\includegraphics[scale=0.1]{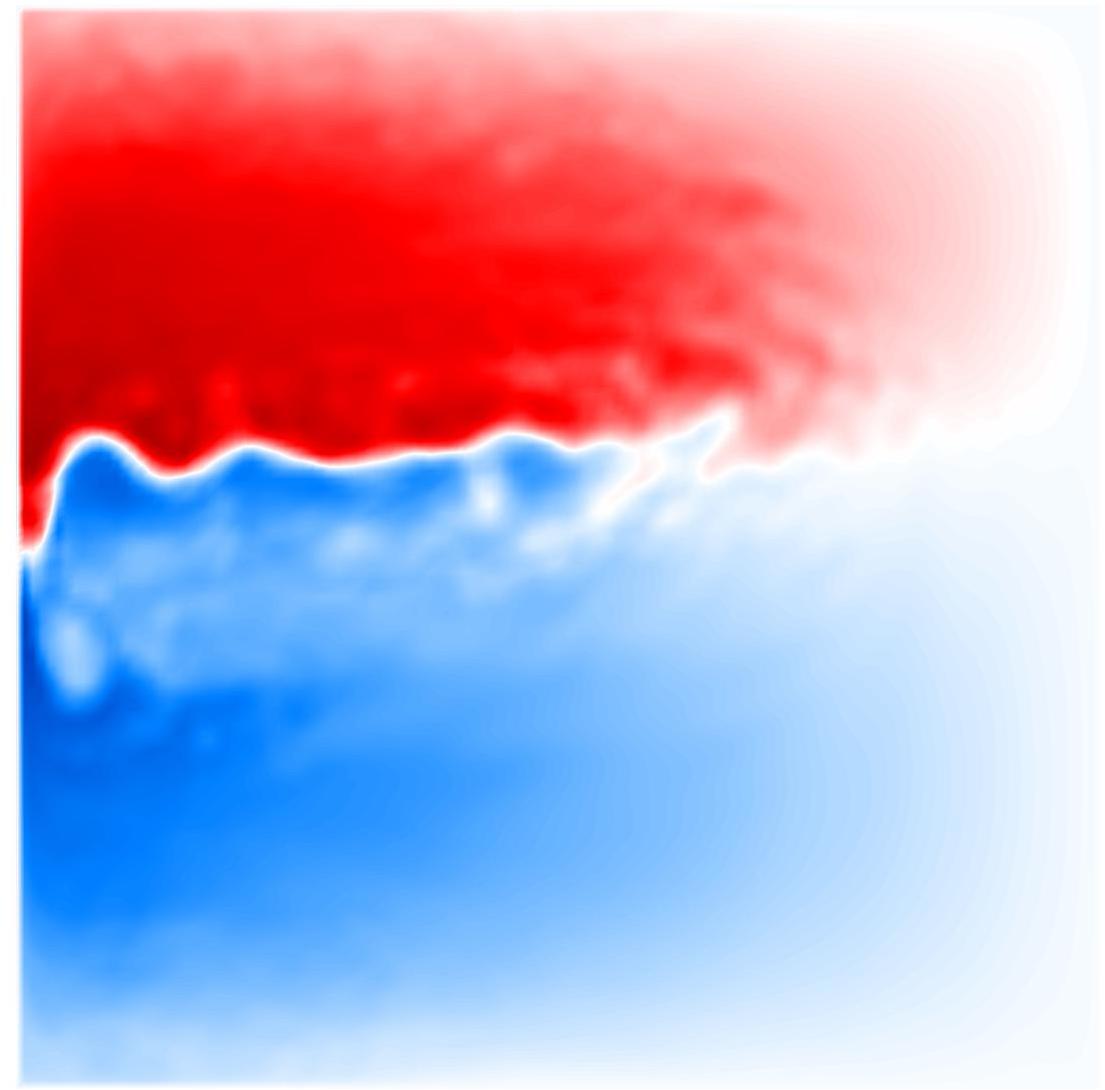}\end{minipage}\\
& & & & \\[-0.35cm]
\multicolumn{5}{c}{\hspace*{-0.5cm}\includegraphics[width=6cm,height=0.5cm]{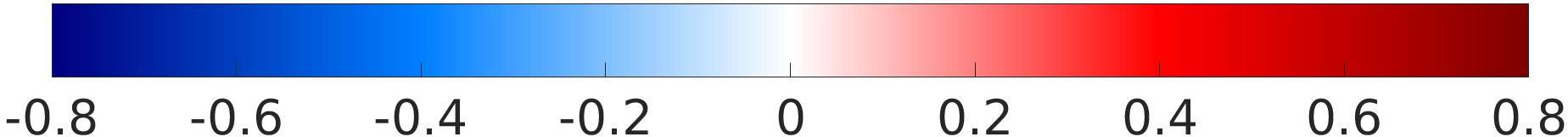}}\\
& & & & \\[-0.15cm]
\hspace*{-1.5cm}\begin{minipage}{0.02\textwidth}\rotatebox{90}{standard deviation}\end{minipage}  &
\hspace*{-1cm}\begin{minipage}{0.24\textwidth}\includegraphics[scale=0.1]{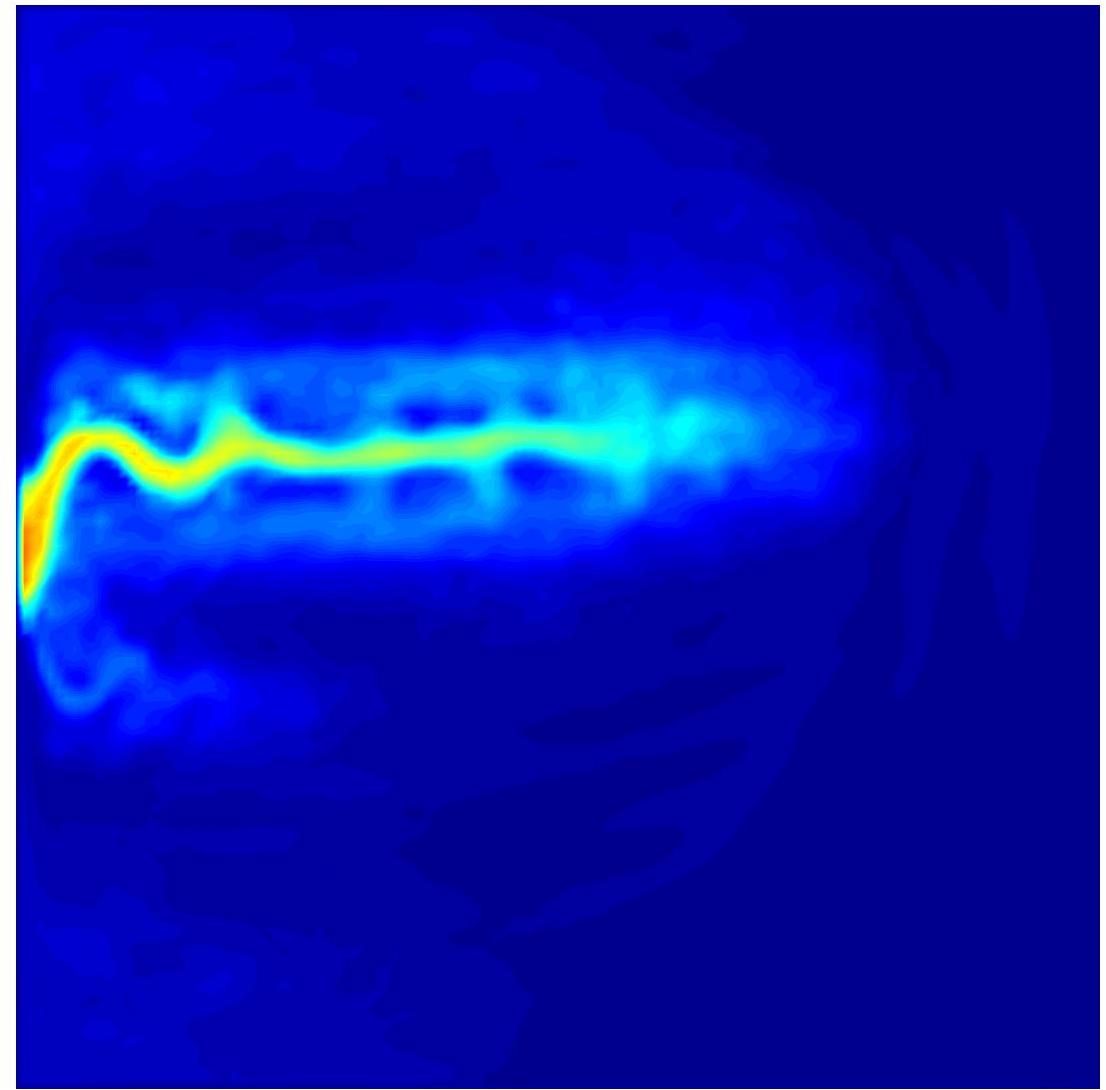}\end{minipage} &
\hspace*{0cm}\begin{minipage}{0.24\textwidth}\includegraphics[scale=0.1]{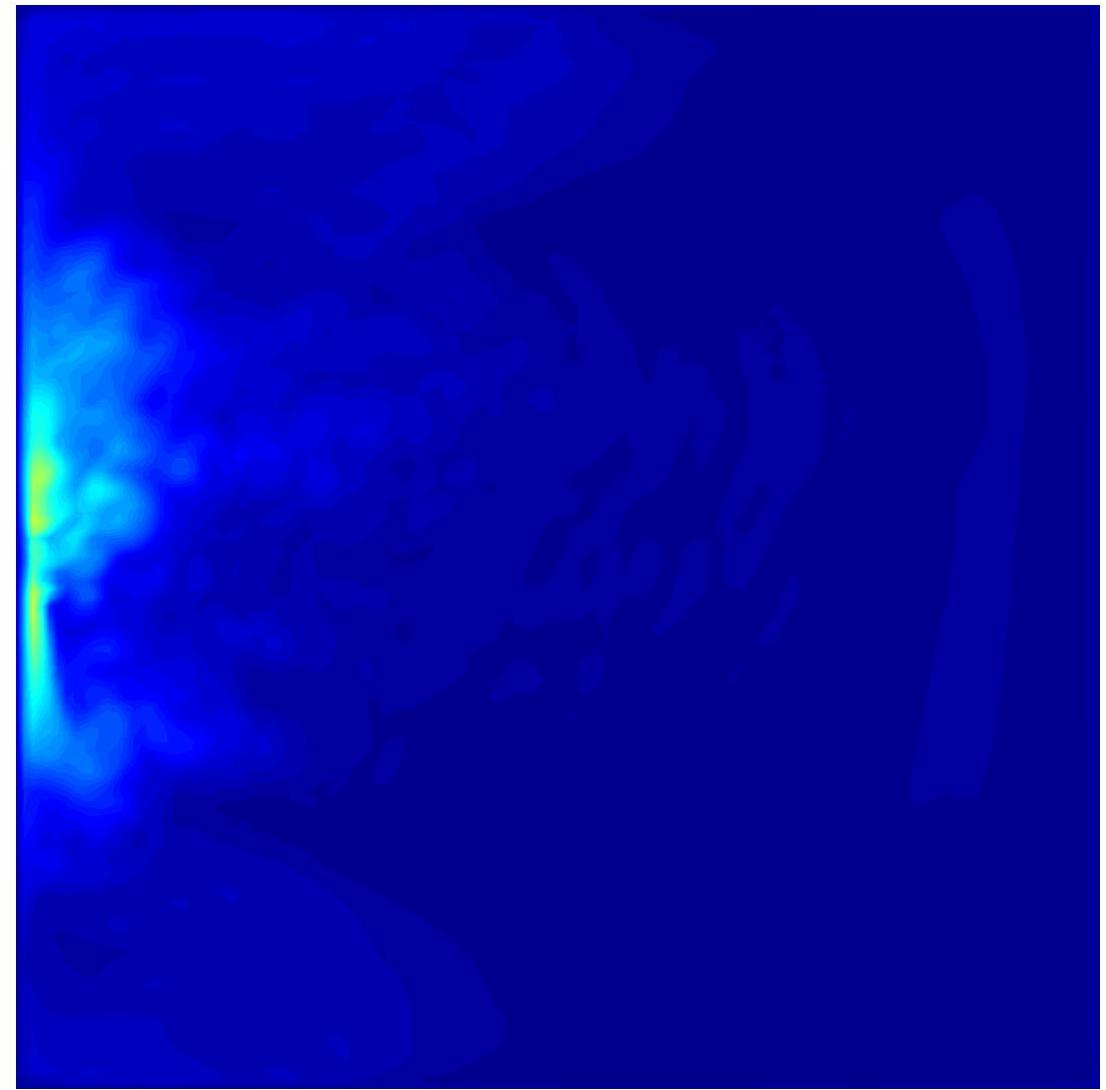}\end{minipage} &
\hspace*{0cm}\begin{minipage}{0.24\textwidth}\includegraphics[scale=0.1]{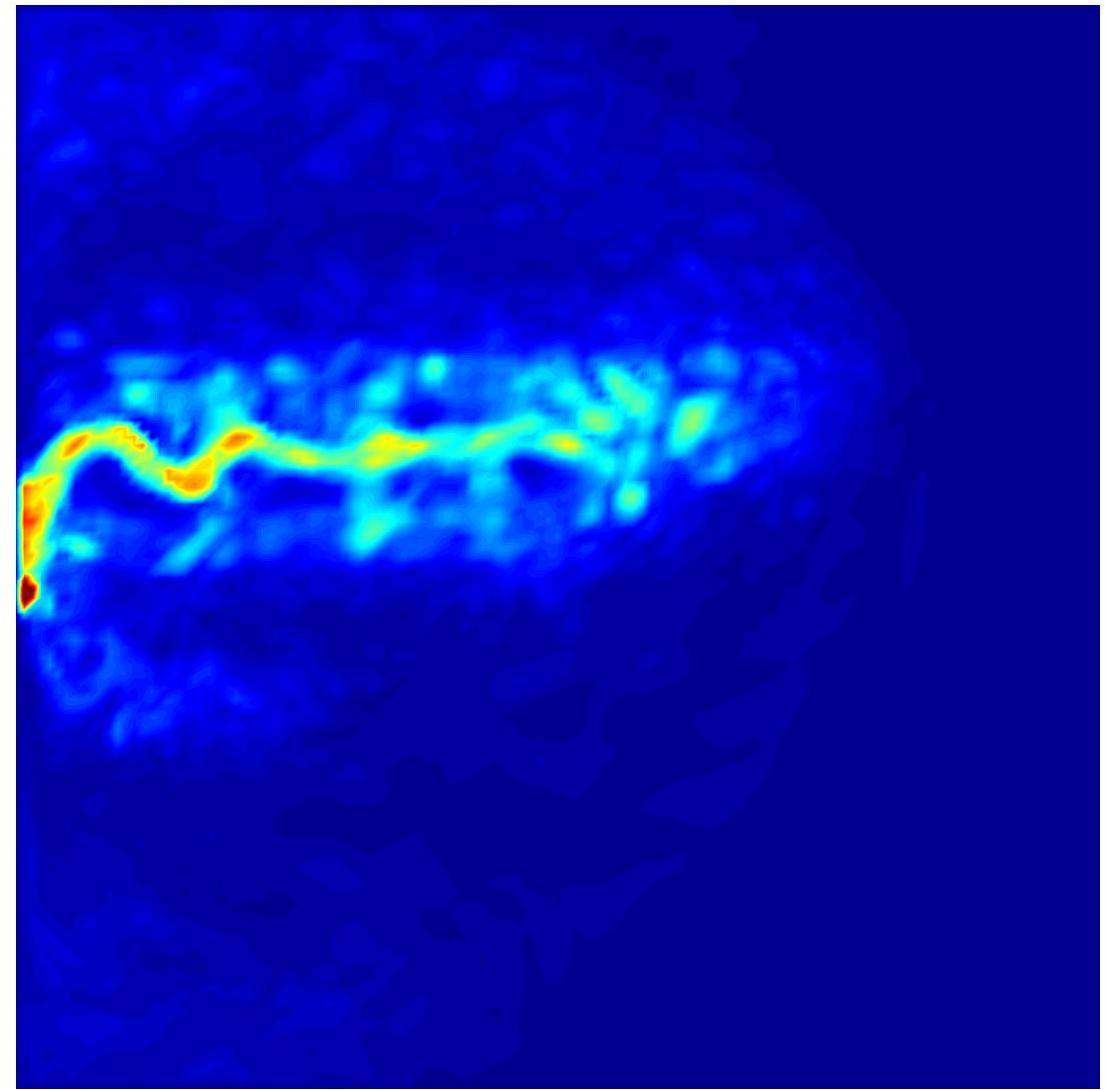}\end{minipage} &
\hspace*{0cm}\begin{minipage}{0.24\textwidth}\includegraphics[scale=0.1]{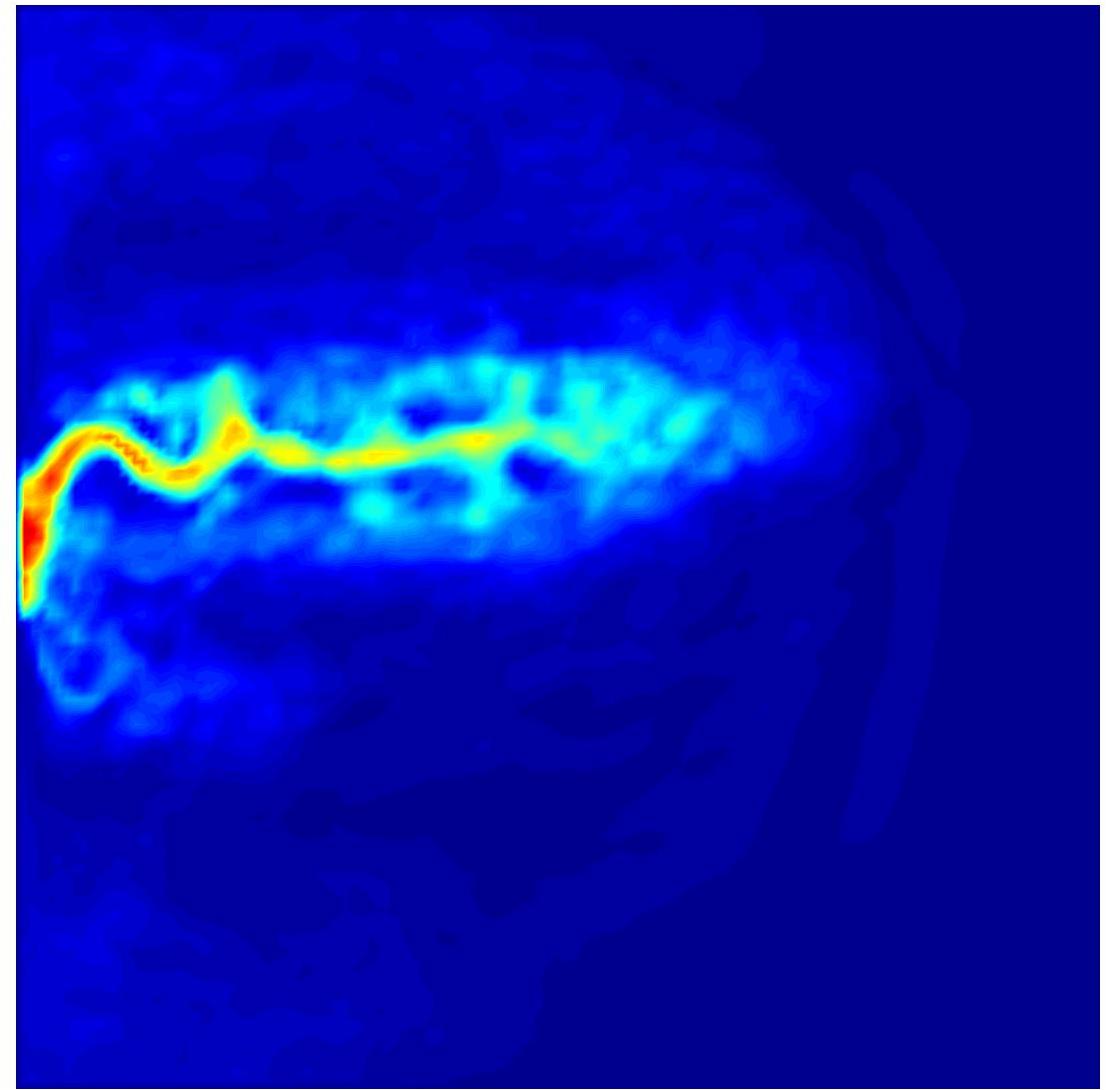}\end{minipage}\\
& & & & \\[-0.35cm]
\multicolumn{5}{c}{\hspace*{-0.5cm}\includegraphics[width=6cm,height=0.5cm]{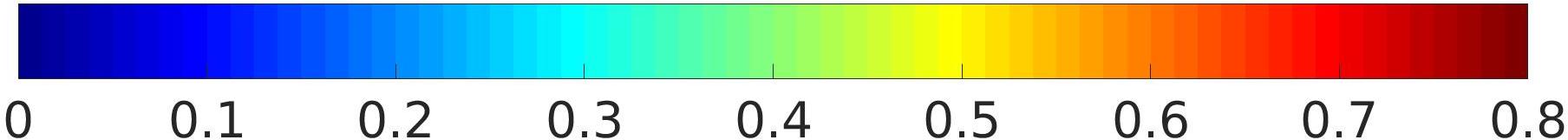}}\\
\end{tabular}
\caption{
Shown is a series of snapshots, 4-year average, and standard deviation of the top layer PV anomaly of {\bf (a)} the reference solution $q_1$ (computed on grid $513\times513$ and projected on grid $129\times129$), {\bf (b)} low-resolution solution 
$\widehat{q}_1$ computed on grid $129\times129$, {\bf (c)} low-resolution solution $\widetilde{q}_1$ on grid $129\times129$ (with the second-order polynomial basis used for the reconstruction), 
{\bf (d)} low-resolution solution $\widetilde{q}_1$ on grid $129\times129$ (with the second-order polynomials and Fourier basis used for the reconstruction).
The solution is given in units of $[s^{-1}f^{-1}_0]$, where $f_0=0.83\times10^{-4}\, {\rm s^{-1}}$ is the Coriolis parameter.
The results in panels {\bf (c)}, {\bf (d)} demonstrate that the proposed method preserves not only large-, but also small-scale features (nominally resolved on the coarse-grid) 
like those seen in the reference solution {\bf (a)} but absent in the low-resolution solution {\bf (b)}.
}
\label{fig:qg_sol}
\end{figure}

For this study we need both high- and low-resolution solutions.
In order to compute them, we first spin up the model~\eqref{eq:pve}-\eqref{eq:bc} for 100 years and then solve it for the other 4 years on 2 uniform horizontal grids: $513\times513$ (high resolution) and $129\times129$ (low resolution).
In order to obtain the reference solution (denoted as $q_1$), we project the high-resolution solution on the coarse grid $129\times129$ by using point-to-point projection (Figure~\ref{fig:qg_sol}a).
The low-resolution solution (denoted as $\widehat{q}_1$) is the solution of the QG model on grid $129\times129$ (Figure~\ref{fig:qg_sol}b).
Our goal is to find a dynamical system that can model the leading PCs (which are then used
to approximate the reference solution given by~\eqref{eq:eof_pc}), so that the approximate solution (denoted as $\widetilde{q}_1$) simulates the reference large-scale flow patterns in qualitatively correct way.

For the purpose of this work, it is enough to consider only the first layer, as it consists of both large- and small-scale features 
(Figure~\ref{fig:qg_sol}a) which we aim to reproduce. Moreover, the upper layer is more difficult to model than the deep ones.
As seen in Figure~\ref{fig:qg_sol}a, the solution is characterized by the well-pronounced eastward jet extension of the western boundary currents and 
surrounding small-scale coherent vortices. 
Both of these features are missed in the low-resolution solution (Figure~\ref{fig:qg_sol}b) due to the under-resolved eddy effects.
In order to restore nominally-resolved flow patterns (the eastward jet and surrounding vortices), 
we first reconstruct a reduced dynamical system (for the leading PCs) which is based on the second-order polynomials and then the one based on the 
second-order polynomials and Fourier series.
\ansA{The solution corresponding to the former} is presented in Figure~\ref{fig:qg_sol}c.
Although the snapshots show that both the eastward jet and vortices are successfully reproduced, 
the time-mean flow significantly differs from the reference solution: 
the eastward jet separation point is shifted north and the jet itself manifests fluctuations unseen in the reference solution.

\begin{figure}[H]
\centering
\hspace*{-4cm}
\begin{tabular}{c}
\hspace*{-0.6cm}\begin{minipage}{0.24\textwidth}\includegraphics[scale=0.1]{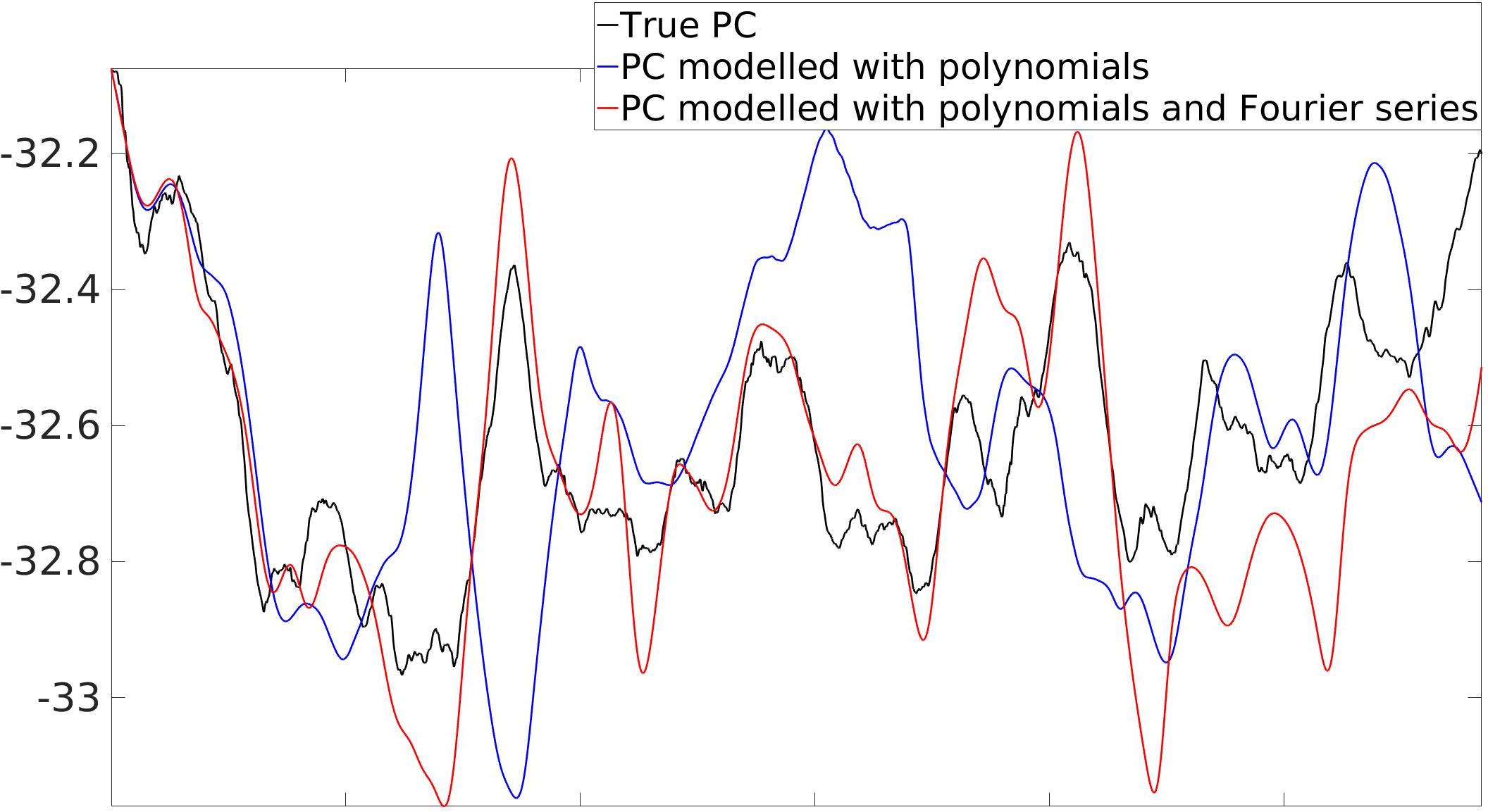}\end{minipage}\\
\\[-0.35cm]
\hspace*{0cm}\begin{minipage}{0.24\textwidth}\includegraphics[scale=0.1]{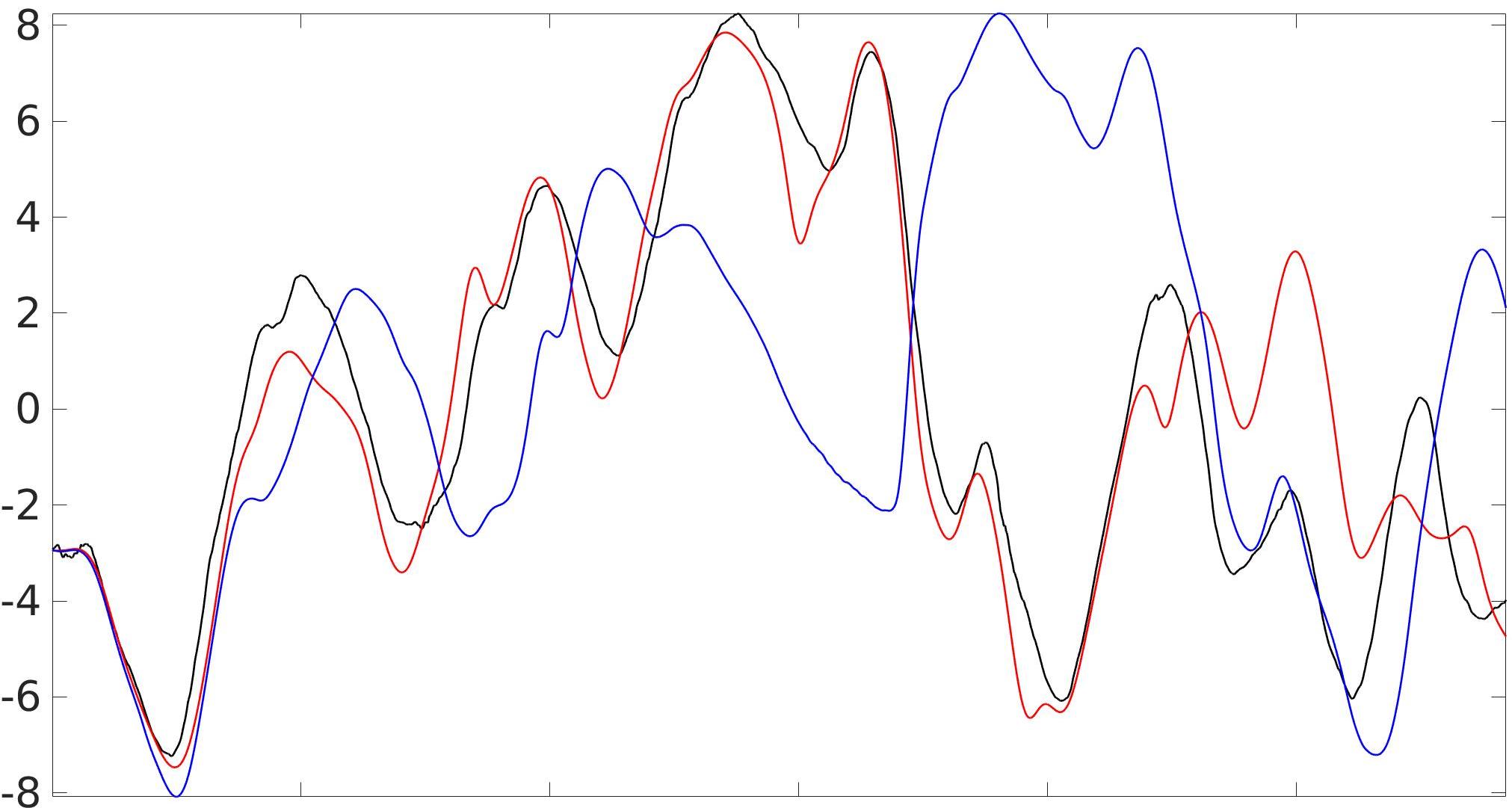}\end{minipage}\\
\\[-0.35cm]
\hspace*{0cm}\begin{minipage}{0.24\textwidth}\includegraphics[scale=0.1]{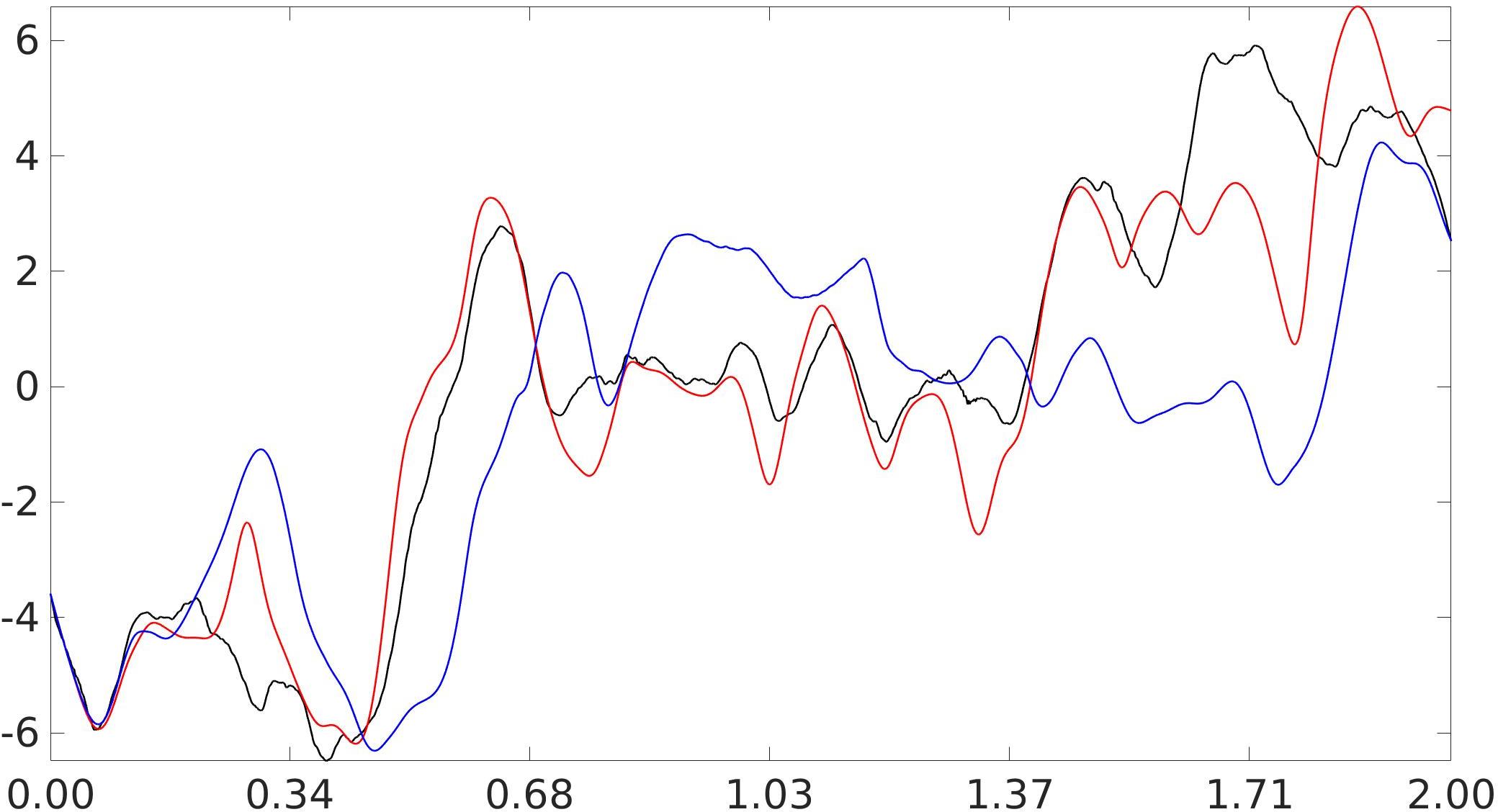}\end{minipage}\\
\\[-0.35cm]
\hspace*{6.75cm}\begin{minipage}{0.24\textwidth}{\bf years}\end{minipage}\\
\end{tabular}
\caption{
Shown are the first three leading PCs and their dependence on the basis functions used for the reconstruction of the dynamical system: true PC (black), PC modelled with the second-order polynomial-only basis (blue), and PC modelled with both the second-order polynomials and Fourier series (red).
The results demonstrate that using the basis consisting of both the second-order polynomials and Fourier series yields significantly more accurate approximation of the PCs.
}
\label{fig:pc}
\end{figure}

For a better approximation we combined the second-order polynomial basis with the Fourier series.
The corresponding solution~\eqref{eq:eof_pc} computed from the leading EOF-PC pairs is significantly improved (Figure~\ref{fig:qg_sol}d),
mostly due to the better approximation of the PCs (Figure~\ref{fig:pc}).

Recall that the solution in Figure~\ref{fig:qg_sol} is over 4 years, and only the first 2 years were used to reconstruct the dynamical system.
This shows that the proposed method preserves not only the large-scale flow structure but also the small-scale flow features, all of them over a long time interval.
The ability of the method to reproduce small-scale features may look surprising, but since these features were present in the reference data, 
their reconstruction is a matter of the high-quality reconstruction of the dynamical system.

A key ingredient that makes the method work is the adaptive nudging which keeps the solution in the right region of the phase space that is occupied by the reference solution.
As an approximation of the reference region, we used a sphere centered at the time-mean of the solution, 
and the sphere radius is the mean distance of the solution from the centre.
The mean distances for the reference and low-resolution solutions are
$\langle\mathcal{D}(q_1,\overline{q}_1)\rangle=11.9$ and 
$\langle\mathcal{D}(\widehat{q}_1,\overline{\widehat{q}}_1)\rangle=7.2$, respectively, showing that the latter is
confined in a smaller region.
The $l_2$-norm distance between the time means of these solutions (denoted as barred quantities)
is $\mathcal{D}(\overline{q}_1,\overline{\widehat{q}}_1)=12.92$.
The application of the adaptive nudging decreases the distance between the time means to 
$\mathcal{D}(\overline{q}_1,\overline{\widetilde{q}}_1)=2.65$, thus shifting the whole solution $\widetilde{q}_1$ much closer to the phase space region occupied by the reference solution.
It also yields a lot more accurate mean distance
$\langle\mathcal{D}(\widetilde{q}_1,\overline{\widetilde{q}}_1)\rangle=12.6$, thus suggesting that the solution has correct amplitude.

\section{Conclusions and discussion\label{sec:conclusions}}
In this study we proposed a method for preserving nominally-resolved flow patterns in low-resolution ocean model simulations. 
The method utilizes the well-known idea of reconstructing the dynamical system that underlies the observed flow evolution.
However, direct application of this idea to the quasi-geostrophic model studied in this work is numerically unfeasible task because of the high dimensionality of the observed flow.
Moreover, a numerical integration of the reconstructed dynamical system can be unstable, but our methodology can cope with this and ensure stability.
We solved the problem of large dimensionality by applying the Empirical Orthogonal Function decomposition of the reference solution
(the high-resolution solution subsampled on the coarse grid) that allowed to reduce the dimension by three orders of magnitude. 
In order to solve the unstable integration problem, we developed the adaptive nudging method following \citep{ShevchenkoBerloff_2021}.
This method keeps the solution in the neighbourhood of the phase space region occupied by the reference solution. 
This is sufficient for accurate reproduction of
both the large- and small-scale flow features at low resolutions, despite the fact that these features are not present in the dynamical solutions of the low-resolution model.
\ansA{The proposed method aims to operate with hundreds of degrees of freedom thus offering orders-of-magnitude acceleration compared 
to low-resolution ocean models which have at least 3-4 orders of magnitude more.
}

The proposed method was tested on a 3-layer quasi-geostrophic ocean circulation model at low non-eddy-resolving resolution, such that it cannot simulate the correct large-scale flow structure.
Our results show that if the reconstructed dynamical model is based only on the second-order polynomials, then it is not sufficiently accurate, because its time-mean eastward 
jet separation point is shifted north, and the jet itself has unrealistic fluctuations which are not observed in the reference solution.
\ansA{We tried to use higher 
order polynomials, but the reconstructed system became very sensitive to errors leading to sever numerical instabilities which we failed to stabilize.}
We resolved this problem by augmenting the polynomial basis with the additional Fourier series.
With all this in place, not only the large-scale flow structure becomes correct but also the small-scale coherent vortices, 
which are unresolved in the low-resolution full-dynamics model, appear in the solution.
All in all, this shows that the method has potential for modelling even more complicated oceanic flows.
Being small-scales-unaware \ansA{(not relying on reproducing the effect of small scales onto large ones like parameterisations)}, 
the proposed method can be thought of as an alternative to the modern (small-scales-aware) parameterisations, which
try to reproduce effects of small dynamically unresolved scales on the large scales, in the hope that the solution will stay in the right region of the phase space.
The proposed approach is quite the opposite: it gently forces the solution to stay in the right phase space region and predicts the flow evolution via the reconstructed reduced dynamical system. 
Note that the method does not require the original quasi-geostrophic model to be solved at low-resolution.

The reference data is used twice: first, for reconstructing the dynamical system; second, for augmenting the solution of this system by nudging.
The method can be further improved by using a more sophisticated equation-wise nudging methodology 
and different dynamical systems which can better represent the underlying flow dynamics.
The proposed method can be straightforwardly applied to primitive equations, but in this case reconstruction of the dynamical system
will be more subtle, as it will include more PCs and can require changes of the basis functions.
Besides, the adaptive nudging may also require some changes, since the phase space behavior of solutions is expected to be more complicated.

\ansA{
Another future extensions of this study can be (1) the requirement to always have a high-resolution simulation (reference solution) 
from which to derive reduced order models -- therefore the ability to compare various existing reference solutions, including observational datasets 
in this framework and to use them systematically/in combination, perhaps without necessarily running fresh high-resolution simulations, 
in order to generate reduced order models at arbitrary resolution;  (2) 
exploring the possibility of generating a forcing for the low-resolution ocean model, based on the EOFs and the adaptive nudging, 
perhaps with some ingredients such as stochastic forcing.
}

\section{Acknowledgments}
The authors thank The Leverhulme Trust for the support of this work through the grant RPG-2019-024.
Pavel Berloff was supported by the NERC grants NE/R011567/1 and NE/T002220/1, and by the Moscow Center of Fundamental and Applied Mathematics (supported by the Agreement 075-15-2019-1624 with the Ministry of Education and Science of the Russian Federation).

\bibliographystyle{apalike}
\bibliography{refs}

\end{document}